\documentclass[final, 1]{elsarticle}
\usepackage{graphicx} 
\usepackage{amssymb, amsmath, amsthm} 
\usepackage{epstopdf}
\usepackage[table]{xcolor}
\usepackage{colortbl}
\usepackage{float}
\usepackage{subfig}
\usepackage{enumerate}

\newtheoremstyle{plainsl}
	{\topsep}
	{\topsep}
	{\slshape} 
	{}
	{\normalfont\bfseries}
	{.}
	{ }
	{}

\renewcommand\proof{\noindent\textsl{Proof. }}
\newcommand\sqr[2]{{\vbox{\hrule height.#2pt
    \hbox{\vrule width.#2pt height#1pt \kern#1pt
        \vrule width.#2pt}\hrule height.#2pt}}}
\renewcommand\qed{%
	\ifmmode\eqno\sqr53
	\else\nolinebreak\ \hfill\sqr53\medbreak\fi}

\newcommand{\col}{\mathrm{col}}
\newcommand{\mur}{\mathrm{mur}}
\newcommand{\mr}{\mathrm{mr}}
\newcommand{\rank}{\mathrm{rank}}

\theoremstyle{plainsl}
\theoremstyle{definition}

\newtheorem{theorem}{Theorem}[section]
\newtheorem{corollary}[theorem]{Corollary}
\newtheorem{lemma}[theorem]{Lemma}
\newtheorem{proposition}[theorem]{Proposition}
\newtheorem{example}[theorem]{Example}

\title{The minimum rank of universal adjacency matrices} 
\begin{document}

\begin{frontmatter}

\author{%
B.~Ahmadi$^{\rm a,c}$
\sep
F.~Alinaghipour$^{\rm a,d}$
\sep
Shaun M. Fallat$^{\rm a,e}$
\sep
Yi-Zheng Fan$^{\rm b}$
\sep \newline
K. Meagher$^{\rm a,f}$
\sep 
S. Nasserasr$^{\rm a,g}$ 
}
 
 \address{$\rm ^a$ Department of Mathematics and Statistics, University of
Regina, Regina S4S0A2, Canada}
\address{$\rm ^b$ School of Mathematical Sciences, Anhui University, Hefei
230039, P.R. China, Email:fanyz@ahu.edu.cn}
 \address{Email: $\rm ^c$~ahmadi2b@uregina.ca, $\rm ^d$~alinaghf@uregina.ca, $\rm ^e$~sfallat@math.uregina.ca, $\rm ^f$~karen.meagher@uregina.ca, $\rm ^g$~shahla.nasserasr@uregina.ca}
\date{}

\begin{abstract}
In this paper we introduce a new parameter for a graph called the
  {\it minimum universal rank}.  This parameter is similar to the minimum rank of a graph. For a graph $G$ the
  minimum universal rank of $G$ is the minimum rank over all matrices of the form 
\[
U(\alpha, \beta, \gamma, \delta) = \alpha A + \beta I + \gamma J + \delta D
\]
where $A$ is the adjacency matrix of $G$, $J$ is the all ones matrix
and $D$ is the matrix with the degrees of the vertices in the main
diagonal, and $\alpha\neq 0, \beta, \gamma, \delta$ are scalars.  Bounds for general graphs based on known graph parameters
are given, as is a formula for the minimum universal rank for regular
graphs based on the multiplicity of the eigenvalues of $A$. The exact value
of the minimum universal rank of some families of graphs are
determined, including complete graphs, complete bipartite graph, paths
and cycles. Bounds on the
minimum universal rank of a graph obtained by deleting a single vertex
are established. It is shown that the minimum universal rank is not
monotone on induced subgraphs, but bounds based on certain induced subgraphs,
including bounds on the union of two graphs, are given. Finally we characterize all graphs with minimum universal rank equal to
$0$ and to $1$.
\end{abstract}

\begin{keyword}
 adjacency matrix \sep universal adjacency matrix \sep Laplacian matrix \sep
minimum rank \sep graph \sep path \sep cycle.\\ 
 \MSC 05C50 \sep 15A03 \sep 15A18 \sep 15A27

\end{keyword}

\end{frontmatter}

\section{Introduction} \label{Sec:intro}

The minimum rank problem for a given graph is a well-studied problem
in the spectral theory of graphs.  The {\em minimum rank of a graph
  $G$} is the smallest rank among all real-valued, symmetric matrices 
that have the property: for $i\neq j$ the $(i,j)$-th entry is
nonzero if and only if $\{i,j\}$ is an edge in the graph $G$. Such a
quantity is denoted by $\mr(G)$. This number has been at the root of a
number of studies over the past dozen years, and a complete resolution
of determining $\mr(G)$ for all $G$ seems essentially unattainable~\cite{FH}.

Haemers and Omidi in~\cite{HO} defined a new family of matrices that
is associated to a graph, these matrices are called the
{\it universal adjacency matrices} of the graph. Consider a simple undirected
graph $G=(V,E)$ with $V=\{v_1,v_2,\ldots, v_n\}$. Let $A_G=[a_{ij}]$
be the $(0,1)$-adjacency matrix of $G$, that is $a_{ij}$ equals one if
$\{i,j\}\in E$ and zero otherwise. Let $D_G=\mathrm{diag}[d_1, d_2,
\ldots, d_n]$ with $d_i=\deg(v_i)$ be the degree matrix associated with
$G$, and denote by $I$ and $J$ the $n\times n$ identity matrix and
$n\times n$ matrix of all ones. An $n\times n$ matrix of the form 
\[
U_G = U_G(\alpha, \beta,\gamma,\delta)=\alpha A_G+ \beta I +\gamma J +\delta D_G,
\]
where $\alpha, \beta, \gamma, \delta$ are scalars with
$\alpha\neq 0$ is called a {\it universal adjacency matrix} of $G$. We
drop the subscript $G$ when it is clear from the context.  The entries
of a universal adjacency matrix $U=[u_{ij}]$ are then of the following
form:
\[
u_{ij} =\left \{
  \begin{array}{cccc} 
  & \beta+\gamma+\delta d_i  &   {\mathrm{if\ } i=j}  \\
 & \alpha+\gamma  &  {\mathrm{if\ } \{i,j\}\in E} \\
  & \gamma  &  {\mathrm{if\ } \{i,j\}\notin E.}
\end{array}
 \right.
\]

Throughout this paper, graphs are considered to be simple and
undirected. Thus a universal adjacency matrix is always a symmetric matrix.
 
The family of universal adjacency matrices is a generalization of several
families of matrices associated to the graph. The following table
shows that for specific values of the coefficients, the universal
adjacency matrix is a well-known matrix associated with a graph:

\renewcommand{\arraystretch}{1.25}
\begin{figure}[h] 
\begin{tabular}{|l|l|} \hline
$(\alpha, \beta, \gamma, \delta)$ &  resulting matrix \\ \hline
$(1,0,0,0)$ & adjacency matrix\\
$(-1,0,0,1)$ & Laplacian matrix\\
$(1,0,0,1)$ & signless Laplacian matrix\\
$(-2, -1, 1,0)$ & Seidel matrix\\
$(\alpha, \beta, \gamma,0)$ &  generalized adjacency matrix\\
$(-1, -1,1,0)$ &  adjacency matrix of the complement \\ \hline
 \end{tabular}
\caption{ Specific cases of the universal adjacency matrix of a graph \label{equivalents}}
\end{figure}

Haemers and Omidi in~\cite{HO} studied the number of distinct
eigenvalues of a universal adjacency matrix associated to $G$ and
determined exactly which graphs have two distinct eigenvalues. In this
paper we are concerned with the rank of universal adjacency matrices.

For a given graph $G$, the {\it minimum universal rank} of $G$, denoted by $\mur(G)$, is given by
\[
\mur (G)=\min \{\rank(U)\ |\ U\ \mathrm{ is\ a\  universal\ adjacency\ matrix\ of}\ G \}.
\]
It is clear that the minimum rank of any of the matrices in
Figure~\ref{equivalents} is an upper bound on the minimum universal
rank of a graph.

If $U_G=\alpha A_G+ \beta I +\gamma J +\delta D_G$ is a universal
matrix for a graph $G$, then, since $\alpha \neq 0$,
\[
A_G+ \frac{\beta}{\alpha} I +\frac{\gamma}{\alpha} J +\frac{\delta}{\alpha} D_G
\] 
is also a universal adjacency matrix for $G$ with the same rank. Thus, in
studying the minimum universal rank of a graph $G$, we may assume
without loss of generality that $\alpha =1$ for a universal adjacency
matrix of $G$.

When considering a universal adjacency matrix, the off-diagonal
entries also come in two types: if an entry corresponds to an edge it
is a fixed number, otherwise it is required to be a different
fixed value. This is in contrast to the off-diagonal entries of the
matrices associated with the minimum rank of $G$, in this case the
entries that correspond to non-adjacent vertices must all be zero,
while the entries corresponding to adjacent vertices are non-zero but
otherwise independent (excluding their symmetric mate). Further, the
main diagonal of a universal adjacency matrix is not completely free
as it is in the matrices associated to minimum rank, but rather it
depends on the degree of a vertex, and the parameters $\beta$,
$\gamma$ and $\delta$. Consequently, the parameters $\mr(G)$ and
$\mur(G)$ are not comparable in general (see examples throughout this
work) and appear to not share any sort of strong relationship. For instance, a graph $G$ in the assumption of Theorem \ref{KrcupKsveevspecial} satisfies $\mathrm{mr}(G)<\mur(G)$, while using Theorem \ref{mur=0}, we have $\mur(K_n)<\mathrm{mr}(K_n).$
However, note that for a given graph $G$, the universal matrix $U_G(1,\beta,\delta,0)$ represents a zero-nonzero pattern for $G$. So if $\mur(G)=\mathrm{rank}(U_G(1,\beta,\delta,\gamma))$, then 
\[
\mathrm{mr}(G)\leq \mathrm{rank}(U_G(1,\beta,\delta,0))=\mathrm{rank}\left (U_G(1,\beta,\delta,\gamma)-\gamma J\right ) \leq \mur(G)+1.
\]

In our notation, $J_{r,s}$ denotes the $r \times s$ matrix of all entries equal to one and $0_{r,s}$ denotes the $r \times
s$ zero matrix. We use $e$ to denote the all ones vector and add
a subscript if it is necessary to specify the size of the vector.

\section{Basic Results}\label{sec:prelim}

In this section we give some basic results about minimum universal
rank for general graphs. The first result shows that the minimum
universal rank has an unusual property that neither the minimum rank
nor the minimum rank of the generalized adjacency matrix has,
namely that the minimum universal rank of a graph is equal to the
minimum universal rank of its complement. We use $\overline{G}$
to denote the complement of the graph $G$.

\begin{lemma}\label{lem:complements}
For any graph $G$, $\mur(G)=\mur(\overline{G})$.
\end{lemma}
\proof  
For any universal adjacency matrix
\[
U_G=\alpha A_G+ \beta I +\gamma J +\delta D_G,
\]
using the facts that $A_{\overline{G}}=-A_G-I-J$ and $D_{\overline{G}}=(n-1)I-D_G$, it follows that
\begin{align*}
U_G &= (-\alpha) (-A_G\!-\!I\!-\!J)+ (\beta\!-\!\alpha\!+\!(n\!-\!1)\delta) I 
        +(\gamma+\alpha) J + (-\delta) ((n\!-\!1)I\!-\!D_G) \\
&= \alpha'(A_{\overline{G}}) + \beta' I + \gamma' J + \delta' D_{\overline{G}}
\end{align*} 
is a universal adjacency matrix for $\overline{G}$.  This implies that for any
set of scalars $\alpha, \beta, \gamma, \delta$ there is a set of
scalars $\alpha',\beta',\gamma',\delta'$ such that the universal
matrices $U_G=\alpha A_G+ \beta I +\gamma J +\delta D_G$ and
$U_{\overline{G}}=\alpha' A_{\overline{G}}+ \beta' I +\gamma' J
+\delta' D_{\overline{G}}$ are equal. Therefore, $G$ and
$\overline{G}$ have the same minimum rank.\qed

The proof of above lemma also shows that the set of universal
adjacency matrices for a graph is equal to the set of universal
adjacency matrices for its complement.

If a graph is disconnected, then its complement is connected, thus in
discussing the minimum universal rank of a graph, we may assume that
the graph is connected, although it may not always be convenient to do so.

The next result shows that it is possible for the minimum universal
rank of a graph to be zero, but this can only happen in a specific case.

\begin{theorem}\label{mur=0}
For any graph $G$, $\mur(G)=0$ if and only if $G$ or $\overline{G}$ is a complete graph.
\end{theorem}
\proof
 Let $\alpha=1$, $\beta =1, \delta=0$ and $\gamma = -1$, then the resulting
  universal adjacency matrix for $K_n$ is the zero matrix.  For the converse, it
  suffices to note that if $G$ has edges and non-edges at the same
  time, then any universal adjacency matrix of $G$ will have a non-zero entry;
  therefore, if $\mur(G)=0$, then $G$ is either complete graph or
  empty graph.\qed

For a graph $G$ on $n$ vertices, the matrix $L_G=D_G-A_G$ is the
Laplacian matrix of $G$.  The Laplacian matrix of a graph is a
universal adjacency matrix of the graph, so the rank of the Laplacian is an
upper bound on the minimum universal rank of the graph. Much is known
about the eigenvalues of a Laplacian matrix that can be used to bound
the minimum universal adjacency matrix of a graph; see \cite{M} for more details.

For example, it is known that $L_G$ is positive
semi-definite. Moreover, the sum of the entries in each row of $L_G$
is zero which implies that zero is an eigenvalue for $L_G$ and $e$, is a corresponding
eigenvector. Furthermore, the multiplicity of zero as an eigenvalue of the Laplacian matrix is
exactly the number of components of the graph. That is, if $c(G)$
denotes the number of components of $G$, and $m_A(\lambda)$ denotes
the multiplicity of $\lambda$ as an eigenvalue of $A$, then we have
$m_{L_G}(0)=c(G)$.

\begin{theorem} \label{upbd-comp}
For any graph $G$ on $n$ vertices
$$ \mur(G)\le n-c(G).\qed$$
\end{theorem}

Note that, the upper bound above cannot be improved since equality holds for the empty graph. One interesting fact about Theorem~\ref{upbd-comp} is that it relates the minimum universal rank of a graph to a
well-known graph parameter, but sometimes it is possible to use the eigenvalues
of the Laplacian matrix to get a better bound on the minimum universal rank of a graph.

\begin{theorem} \label{upbd-multi} 
  Let $G$ be a connected graph on $n$ vertices, and let $m$ be the maximum
  multiplicity of the nonzero eigenvalues of $L_G$.  Then
\[
 \mur(G)\le n-m-1.
\]
\end{theorem}

\proof Let the eigenvalues of $L_G$ be $\lambda_1=0\leq \lambda_2\leq 
\ldots\leq \lambda_n$.  The all ones vector
is an eigenvector for $0$ and the eigenvectors for the nonzero
eigenvalues are orthogonal to the all ones vector. Suppose
$\lambda_k$ is a nonzero eigenvalue of $L_G$ with multiplicity $m$,
then $\lambda_k$ is an eigenvalue of $L_G+\frac{\lambda_k}{n}J$ with
multiplicity $m+1$ (to see this, consider the $m$ linearly independent
eigenvectors for $L_G$ corresponding to $\lambda_k$ and the all ones
vector).  So, the matrix $L_G+\frac{\lambda_k}{n}J-\lambda_k I$ is a
universal adjacency matrix for $G$ that has rank $n-m-1$, and hence $\mur(G) \le
n-m-1$. \qed

Note that, if $G$ is connected, then the multiplicity of the eigenvalue zero of the Laplacian matrix is one (the number of components). So if zero has the maximum multiplicity, then all eigenvalues are simple, which implies $m=1$ in above proof. We also note that Theorem~\ref{upbd-multi} is valid even if $G$ is disconnected, since then its complement $\overline{G}$ is connected and has the same minimum universal rank as $G$. So if $G$ is disconnected we apply the above proof for $\overline{G}$.

The following is an immediate consequence of Theorem \ref{upbd-multi}. 

\begin{corollary} \label{upbd}
For any graph $G$ on $n$ vertices
\[
 \mur(G)\le n-2.\qed
\]
\end{corollary}

It is known that the minimum rank of a graph on $n$ vertices is at most $n-1$, and mr$(G)=n-1$ if and only if $G$ is a path on $n$ vertices; see \cite{FH}. In the next section, it is shown that the upper bound in Corollary~\ref{upbd} is achieved by paths. It is interesting to note that the graphs that achieve the
maximum possible minimum rank also achieve the maximum possible
minimum universal rank. But unlike minimum rank, where the paths are
the only graphs that have the maximum possible minimum rank, there are
many graphs that achieve the maximum possible minimum universal rank; see Example \ref{lem:paths} for instance.
 
\section{Paths}

A {\it path} on $n$ vertices, denoted by $P_n$, is a graph with
vertices $v_1,v_2,\ldots,v_n$ and edge set 
$\{\{v_1,v_2\}, \{v_2,v_3\},\ldots, \{v_{n-1},v_n\}\}$.  The next
result shows that if a graph contains an induced path on $n$ vertices, then the
minimum universal rank of the graph is at least $n-2$.

Let $A$ be an $m\times n$ matrix. For $\alpha \subseteq
\{1,2,\ldots,m\}$ and $\beta \subseteq \{1,2,\ldots,n\}$, the notation
$A[\alpha,\beta]$ means the submatrix of $A$ lying in rows indexed by
$\alpha$ and columns indexed by $\beta$.

\begin{lemma}\label{lbd-indpath}
If a graph $G$ contains the induced path $P_k \; (k \ge 3)$, then 
\[
\mur(G) \ge k-2.
\]
\end{lemma}

\proof
Suppose $v_1v_2\ldots v_k$ is the induced path $P_k$.
Order the vertices of $G$ such that the first $k$ vertices are
$v_1,v_2,\ldots,v_k$. Then any universal adjacency matrix $U_G=A_G+\beta I
+ \gamma J + \delta D_G$ of the graph $G$ is of the form:
\[
\left[\begin{array}{c|c}
{
\begin{array}{ccccc}
\beta\!+\!\gamma\!+\!\delta d_{v_1} & 1\!+\!\gamma & \cdots & \gamma & \gamma \\
1\!+\!\gamma & \beta\!+\!\gamma\!+\!\delta d_{v_2}&  \cdots & \gamma & \gamma \\
\gamma & 1\!+\!\gamma &  \cdots  & \gamma & \gamma \\
\vdots & \vdots & \ddots & \vdots & \vdots \\
\gamma & \gamma &  \cdots   & \beta\!+\!\gamma\!+\! \delta d_{v_{k-1}} & 1\!+\!\gamma \\
\gamma & \gamma &  \cdots   & 1\!+\!\gamma & \beta\!+\!\gamma\!+\! \delta d_{v_k} \\
\end{array}
} & B \\
\hline 
B^T & C 
\end{array}\right],
\]
Subtracting the $k$th column from each of the columns  $2, \ldots, k-1$ results in the following matrix 
\[
U'=\left[\begin{array}{c|c}
{
\begin{array}{cccccccc}
\beta\!+\gamma\!+\delta d_{v_1} & 1 & 0 & &  \cdots  & & 0 & \gamma \\
1+\gamma & * & 1 & & \cdots & & 0 & \gamma \\
\gamma & * &* & &  \cdots & & 0 & \gamma \\
\vdots & \vdots & & & \vdots & & \vdots & \vdots \\
\gamma &  & & & \cdots  & & 1  & \gamma \\
\gamma &  & &  &  \cdots  & & * & 1+\gamma \\
\gamma &  &  & & \cdots  &  & * & \beta\!+\gamma\!+ \delta d_{v_k} \\
\end{array}
} & B \\
\hline 
* & C 
\end{array}\right].
\]

Since the submatrix $U'[\{1,\ldots,k-2\},\{2,\ldots,k-1\}]$ of $U$
has rank $k-2$, we have $\mathrm{rank}(U_G)\geq\mathrm{rank}(U') \ge
k-2$, which implies $\mur(G) \ge k-2.$ \qed

The {\it diameter} of a graph $G$, denoted by $\mathrm{diam}(G)$, is
the maximum distance between vertices of the graph. Since a path
corresponding to the diameter is an induced path, we have the
following consequence.

\begin{corollary}
For any graph $G$, 
\[
 \mur(G) \geq \mathrm{diam}(G)-1.\qed
\]
\end{corollary}

The {\it union} of graphs $G_1=(V_{1},E_{1})$, $G_2=(V_{2},E_{2}),
\ldots, G_{m}=(V_{m},E_{m})$ is the graph 
\[
\bigcup_{i=1}^{m} G_{i}= \left(\bigcup_{i=1}^{m} V_{i}, \bigcup_{i=1}^{m} E_{i}\right).
\] 
If $G_1=G_2=\cdots=G_{m}$, then instead of $\cup_{i=1}^{m}
G_{i}$, we use the notation $m G$.

\begin{lemma}\label{lbd-indpath3}
If a graph $G$ on $n$ vertices contains the induced subgraph 
$P_{k_1} \cup P_{k_2} \cup \cdots \cup P_{k_t}$ with $t \ge 2$, and $k_i \ge 2, i=1,\ldots,t$, then 
\[
\mur(G) \ge \left(\sum_{i=1}^t k_i \right) -(t+1).
\]
\end{lemma}
\proof Order the paths as given and in each path order the vertices so that a pendant vertex comes first and every other vertex comes right after its previous neighbour. Then any universal adjacency matrix of $G$ is of the following form

\begin{align*}
U=\left[ \begin{array}{c|c|c|c|c}
U_{11} & \gamma J_{k_1, k_2}  & \ldots & \gamma J_{k_1, k_{t}} & U_{1(t+1)}\\
\hline
\gamma J_{k_2, k_1} & U_{22} &  \ldots & \gamma J_{k_2, k_t} & U_{2(t+1)}\\
\hline
\vdots & \vdots &  \ddots & \vdots & \vdots\\
\hline
\gamma J_{k_t, k_{1}} & \gamma J_{k_t, k_2}  & \ldots & U_{tt} & U_{t(t+1)}\\
\hline
U_{1(t+1)}^T & U_{2(t+1)}^T  & \ldots & U_{t(t+1)}^T & U_{(t+1)(t+1)}
\end{array}
\right]. 
\end{align*}
For each block $U_{ii}$, $i=1,\ldots,t$, the super diagonal entries are $\gamma+1$ and every other non-diagonal entry equals $\gamma$. Now subtracting the column $\displaystyle \sum_{i=1}^t k_i$ from each of the columns $1,2, \ldots, \displaystyle\sum_{i=1}^t k_i-1$, we produce a lower triangular submatrix of size $k_i-1$ in the $(i,i)$ block for $i=1,\ldots, k_t-1$ and a lower triangular submatrix of size $k_t-2$ in the
$(t,t)$ block, with all ones on the main diagonal entries of each of the triangular matrices. Moreover, the entries of the blocks above these triangular matrices in the resulting matrix are all zero except possibly the last column. So
$\mathrm{rank}(U_G) \ge \displaystyle \sum_{i=1}^t k_i -(t+1)$. \qed

\begin{lemma}\label{lbd-indpath2}
If $G$ contains the induced subgraph 
$P_{k_1} \cup P_{k_2} \cup \cdots \cup P_{k_t} \cup m P_1$ with $t \ge 2, m \ge 1$, and $k_i \ge 2,  i=1,\ldots,t$, then 
\[
\mur(G) \ge \left(\sum_{i=1}^t k_i \right) -t.
\]
\end{lemma}

\proof The proof is similar to the proof of Lemma \ref{lbd-indpath3}. Using the same ordering for vertices, and subtracting the column corresponding to the column of one of the vertices in $mP_1$, results in a lower triangular of rank $k_t-1$ matrix for the block corresponding to $P_{k_t}$ as well, which implies the inequality.  \qed

In some cases in Lemmas \ref{lbd-indpath3} and \ref{lbd-indpath2}, the equality holds when the induced subgraph $P_{k_1} \cup P_{k_2} \cup \cdots \cup P_{k_t} \cup m P_1$, $(m\geq 0)$, is exactly the graph $G$. Some of these cases are listed below. 

\begin{example}\label{lem:paths} For $n\ge 3$ 
\begin{enumerate}[(a)]
\item $\mur(P_n)=n-2$; \label{path}
\item $\mur(P_{n-1} \cup P_1) = n-2$; \label{pathwithpoint}
\item $\mur(P_n \cup P_n)    = 2n-3$. \label{samepathtwice}
\end{enumerate}
\end{example}
\proof Equation~(\ref{path}) can be obtained from Lemma~\ref{lbd-indpath}
and Corollary~\ref{upbd}.  Moreover, the universal adjacency matrix $-A-\lambda
I + \frac{\lambda}{n}J+D$ meets this bound, where $\lambda$ is an
arbitrary nonzero eigenvalue of $L_{P_n}$. For instance, choosing
$\lambda=2\left(1- \cos(\frac{\pi}{n})\right)$, we get the following universal adjacency matrix of rank
$n-2$
\[ 
U_{P_{n}} = A_{P_{n}} +2\left(1- \cos(\frac{\pi}{n})\right)I - \frac{2}{n}\left(1- \cos(\frac{\pi}{n})\right)J- D_{P_{n}}.
\]

Equation~(\ref{pathwithpoint}) is obtained from Lemma~\ref{lbd-indpath3} and Corollary~\ref{upbd}.
The Laplacian matrix $L_{P_{n-1} \cup P_1}$ is an example of a universal adjacency matrix for the graph that has the minimum rank.

Finally, Equation~(\ref{samepathtwice}) is obtained from
Theorem~\ref{upbd-multi} and Lemma~\ref{lbd-indpath3}.  The universal
adjacency matrix $-A-\lambda I + \frac{\lambda}{n}J+D$ meets this
bound, where $\lambda$ is an arbitrary nonzero eigenvalue of
$L_{P_n}$.\qed

\section{Regular Graphs}\label{Sec:regular}

A graph $G$ is called {\em regular of degree $r$} if each vertex of
$G$ is adjacent to exactly $r$ vertices.  If $G$ is a regular graph of
degree $r$, then it is evident that $A_{G}e=re$, $A_{G}J=JA_{G}$. Moreover,
$D_{G} =rI$, so any universal adjacency matrix associated with a
regular graph may be reduced to the form $U_{G} = A_{G} + \beta I +
\gamma J$, which is the generalized adjacency matrix of $G$. Now we
are able to derive the following result regarding the minimum
universal rank of any regular graph.

\begin{theorem}\label{maxmultregular}
  Suppose $G$ is a connected $r$-regular graph of degree $r$ on $n$ vertices.
  Let the spectrum of the adjacency matrix, $A_{G}$, be given by $r,
  \lambda_{2}, \lambda_{3}, \ldots, \lambda_{n}$ (these values may not
  be distinct), and assume that $m$ is the maximum multiplicity among
  the list $\{\lambda_{2}, \lambda_{3}, \ldots, \lambda_{n}\}$. Then
\[ \mur(G) = n-(m+1).\]
\end{theorem}

\proof Since $G$ is a regular graph of degree $r$, we know that
$A_{G}e=re$. Thus for any other eigenvalue $\lambda_{i}$ if $x_i$ is a
corresponding eigenvector, then $x_{i}$ is orthogonal to $e$ and hence
$Jx_{i}=0$.  Furthermore, since $A_{G}$ and $J$ commute and are
symmetric, it follows that the eigenvalues of $A_{G} + \gamma J$ are
given by $r+\gamma n, \lambda_{2}, \lambda_{3}, \ldots, \lambda_{n}$.
(Here we also used the facts that $J$ is rank one and $Je=ne$.)  Thus
the maximum number of zero eigenvalues admitted by any universal
matrix $U_{G} = A_{G} + \beta I + \gamma J$ is equal to $m+1$ by
suitable choices of $\beta$ and $\gamma$.  From which it follows that
$\mur(G) = n-(m+1)$, which completes the proof.  \qed

The adjacency eigenvalues of $K_n$ are $\{n-1, -1,-1, \ldots, -1\}$,
where $-1$ occurs with multiplicity $n-1$, in this case $m=n-1$, and
$\mur(K_{n})=0$. This gives an alternative proof of one of the directions in Theorem \ref{mur=0}.

Since the adjacency eigenvalues of a cycle on $n$ vertices are known, we have the following
as an immediate consequence.

\begin{corollary}\label{n-cycle}
For any $n\geq 3, k\geq 1$, 
\[
\mur(kC_{n}) = kn-2k-1. 
\]
\end{corollary}
\proof The adjacency eigenvalues of $C_{n}$ are well known to be twice
the real parts of the $n$-th roots of unity; see \cite{CDS}. Thus, the maximum
multiplicity among the eigenvalues different from the degree is
$m=2$.  Moreover, $kC_n$ has the same eigenvalues as $C_n$ and each
eigenvalue has multiplicity $k$ times its multiplicity as an eigenvalue
for $C_n$. Applying Theorem \ref{maxmultregular}, we have that
$\mur(kC_{n})=kn-(2k+1)$.  \qed

In particular,  this means that $\mur(C_n) = n-3$.

There are many large families of graphs for which all eigenvalues of their adjacency matrices are
known and for these it is easy to determine the minimum universal
rank. For example, the adjacency matrix of any strongly regular graph
has exactly three distinct eigenvalues and the multiplicities of those that are not equal to the degree can be
expressed in terms of the parameters for the graph; see \cite{BR}; thus the minimum
universal rank can also be expressed in terms of the parameters of the
strongly regular graph.

\section{Unions of Graphs}

We have seen several examples of the minimum universal rank for a
graph that is the union of smaller graphs. This motivates us to
consider bounds on the minimum universal rank of the union of two
graphs; we start with what is a natural lower bound. Indeed we are considering unions of graphs as opposed to joins (see definition on page \pageref{join}), as we feel this approach eases exposition.

\begin{lemma}\label{lem:upperboundunion}
Let $G$ and $H$ be two graphs, then
\[
 \mur(G)+\mur(H) \leq \mur(G \cup H).
\]
\end{lemma}
\proof
Suppose $G$ has $n$ vertices and $H$ has $m$ vertices, 
 and suppose that $\mur(G \cup H)$ is attained by the following universal adjacency matrix:
\begin{align}
U_{G \cup H}(1,\beta,\delta,\gamma)=\left[ \begin{array}{c|c}
U_G & \gamma J_{n,m} \\
\hline 
\gamma J_{m,n} & U_H
\end{array}
\right]. \label{eq1}
\end{align}
If $\gamma=0$, then 
$$ \mur(G \cup H) =\rank(U_G)+\rank(U_H)\ge \mur(G)+\mur(H).$$
Thus we may assume $\gamma \ne 0$. We consider two cases, first if the
all ones vector is in the column space of either $U_G$ or $U_H$ and
second when it is not. The column space of a matrix $A$ is denoted by $\col(A)$.

\underline{case 1:} $e_n \in \col(U_G)$ or $e_m \in \col(U_H)$. \\ \hfill
We only consider the case $e_n \in \col(U_G)$, as the other case is similar.
Note that $\gamma e_n$ must be a linear combination of the columns of $U_G$.
Subtracting this combination from each column of $\gamma J_{n,m}$, we arrive at the following matrix, 
where $\zeta \over \gamma$ is the sum of the coefficients of the above linear combination.
\[
\left[ \begin{array}{c|c}
U_G &  0_{n,m} \\
\hline
\gamma J_{m,n} & U_H-\zeta J_m
\end{array}
\right].
\]
As $U_G$ is symmetric, subtracting the corresponding linear combination of the rows of $U_G$ from each 
row of $\gamma J_{m,n}$ we arrive at the matrix:
\[
\left[ \begin{array}{c|c}
U_G &  0_{n,m} \\
\hline
0_{m,n} & U_H-\zeta J_m
\end{array}
\right].
\]
So $$ \mur(G \cup H) =\rank(U_G)+\rank(U_H-\zeta J_m)\ge \mur(G)+\mur(H).$$
Furthermore, if $\zeta \ne 0$ and $e_m \notin \col(U_H)$, then
$\rank(U_H-\zeta J_m)=\rank(U_H)+1$, and hence 
$$ \mur(G \cup H) \ge \mur(G)+\mur(H)+1.$$

\underline{case 2:} $e_n \notin \col(U_G)$ and $e_m \notin \col(U_H)$. \\ \hfill
Applying elementary row operations to the matrix in (\ref{eq1}) on the rows corresponding to $U_G$,
and if necessary, permuting some columns of $U_G$ (still preserving the zero-nonzero pattern of other blocks),
we have the following matrix for some non-singular diagonal matrix 
$\Lambda_{n'}$ of order $n'$, some matrix $B$ of order $n'\times (n-n')$,
and some real numbers $a_1,a_2,\ldots,a_n$.
\[
\left[ \begin{array}{c|c|c}
\Lambda_{n'} & B & \begin{array}{c} a_1 e_m^T \\ \vdots \\ a_{n'} e_m^T \end{array} \\
\hline
0_{n-n',n'} & 0_{n-n',n-n'} & \begin{array}{c} a_{n'+1} e_m^T \\ \vdots \\ a_{n} e_m^T \end{array} \\
\hline
\gamma J_{m,n'} & \gamma J_{m,n-n'} &  U_H
\end{array}
\right].
\]
We find that $n'<n$ as $U_G$ is singular (otherwise $e_n \in \col(U_G)$),  
and at least one  of $a_{n'+1}, \ldots,a_n$, say $a_n$, is nonzero,  as $e_n \notin \col(U_G)$.
Applying row operations again, we get the following matrix:
\begin{align}
\left[ \begin{array}{c|c|c}
\Lambda_{n'} & B & 0_{n',m}\\
\hline
0_{n-n'-1,n'} & 0_{n-n'-1,n-n'} & 0_{n-n'-1,m} \\
\hline
0_{1,n'} & 0_{1,n-n'} & e_m^T \\
\hline
\gamma J_{m,n'} & \gamma J_{m,n-n'} &  U_H
\end{array}
\right]. \label{eq2}
\end{align}

Similarly, by applying some elementary row operations to the
matrix in (\ref{eq2}) on the rows corresponding to $U_H$, and if
necessary, permuting some columns of $U_H$ (still preserving the
zero-nonzero pattern of other blocks), we have the following matrix for some non-singular
diagonal matrix $\bar{\Lambda}_{m'}$ of order $m'<m$, and some matrix
$C$ of size $m' \times (m-m')$.
\begin{align}
\left[ \begin{array}{c|c|c|c}
\Lambda_{n'} & B & 0_{n',m'} & 0_{n', m-m'}\\
\hline
0_{n-n'-1,n'} & 0_{n-n'-1,n-n'} & 0_{n-n'-1,m'} &  0_{n-n'-1,m-m'}\\
\hline
0_{1,n'} & 0_{1,n-n'} &  e_{m'}^T &  e_{m-m'}^T \\
\hline
0_{m',n'} & 0_{m',n-n'} & \bar{\Lambda}_{m'} & C \\
\hline
0_{m-m'-1, n'} & 0_{m-m'-1, n-n'} & 0_{m-m'-1, m'}  & 0_{m-m'-1, m-m'} \\
\hline
 e_{n'}^T &  e_{n-n'}^T & 0_{1,m'} & 0_{1,m-m'} 
\end{array}
\right]. \label{eq3}
\end{align}
By deleting the zero rows from the matrix in (\ref{eq3}), we have
\[
\left[ \begin{array}{c|c|c|c}
\Lambda_{n'} & B & 0_{n',m'} & 0_{n', m-m'}\\
\hline
0_{1,n'} & 0_{1,n-n'} &  e_{m'}^T &  e_{m-m'}^T \\
\hline
0_{m',n'} & 0_{m',n-n'} & \bar{\Lambda}_{m'} & C \\
\hline
 e_{n'}^T &  e_{n-n'}^T & 0_{1,m'} & 0_{1,m-m'} 
\end{array}
\right]. 
\]
So
\belowdisplayskip=-12pt
\begin{align*}
\mur(G \cup H) & =\rank
\left[ \begin{array}{cc} \Lambda_{n'} & B \\ e_{n'}^T &  e_{n-n'}^T \end{array} \right]
+\rank\left[ \begin{array}{cc} e_{m'}^T &  e_{m-m'}^T  \\ \bar{\Lambda}_{m'} & C \end{array} \right]
\\
& = (\rank\Lambda_{n'}+1)+(\rank\bar{\Lambda}_{m'}+1) \\
& \ge \mur(G)+\mur(H)+2. 
 \end{align*} \qed

\belowdisplayskip=12pt
Note that, from the proof we can conclude that if $\gamma e \not \in \col(U_G)$ and
$\gamma e \not \in \col(U_H)$, then we actually have a stronger bound on $\mur(G \cup H)$.

\label{join} The {\it join} of graphs $G_1=(V_{1},E_{1})$, $G_2=(V_{2},E_{2}),
\ldots, G_{m}=(V_{m},E_{m})$ is a graph on the vertices
$\cup_{i=1}^{m} V_{i}$ that includes the edges $\cup_{i=1}^{m} E_{i}$
but also has all edges $\{v_i,v_j\}$ where $v_i \in V_i$ and $v_j \in
V_j$ with $i\neq j$. The join is denoted by $G_1 \vee G_2 \vee
\cdots \vee G_m$. The join and the union are complementary operations
in the sense that for any pair of graphs $G_1$ and $G_2$,
\[
\overline{G_1 \cup G_2} = \overline{G_1} \vee \overline{G_2}.
\]
This fact, together with Lemma~\ref{lem:complements}, yields the following
\[
 \mur(G_1 \cup G_2) = \mur(\overline{G_1 \cup G_2}) = \mur(\overline{G_1} \vee \overline{G_2}).
\]
This means that results about the union of graphs can be translated to results about
joins of graphs, for example the next result is
Lemma~\ref{lem:upperboundunion} stated for the join of two graphs

\begin{lemma}
For graphs $G$ and $H$
\[
\mur(G \vee H) \ge \mur(G)+\mur(H).
\]
\end{lemma}
\proof With the comments above, we simply note that
\[
\mur(G \vee H) = \mur(\overline{G} \cup \overline{H}) \ge \mur(\overline{G}) + (\overline{H}) = \mur(G) + \mur(H). \qed
\]

Upper bounds on the minimum universal rank of the union of two graphs
seem to be more difficult question.  We have seen examples where it is
possible to express the minimum universal rank of the union of graphs
in terms of the minimum universal ranks of the graphs in the
union. For example, (\ref{pathwithpoint}) and
(\ref{samepathtwice}) of Example~\ref{lem:paths} state that
\[
\mur(P_{n-1} \cup P_1) = \mur(P_{n-1})+\mur(P_1) + 1
\]
and 
\[
\mur(P_n \cup P_n)= 2\;\mur(P_n)+1.
\]
From Corollary~\ref{n-cycle}, 
\[
\mur(C_n \cup C_n) = 2\;\mur(C_n)+1,
\] 
and more generally that $\mur(kC_n)= k\mur(C_n)+k-1$.

This might lead one to conjecture that the minimum universal rank of
the union of the two graphs is bounded above by the sum of the minimum universal ranks
of the graphs in the union plus one, but the difference between
$\mur(G \cup H)$ and $\mur(G)+\mur(H)$ can be arbitrarily large. For example, take $G = kC_3$
and $H=kC_4$, so $\mur(kC_3) = k-1$ and $\mur(kC_4) = 2k-1$. But
Theorem~\ref{maxmultregular} implies that $\mur(kC_3 \cup kC_4) = 5k-1$ so
\[
\mur(kC_3 \cup kC_4)  - (\mur(kC_3) +\mur(kC_4) ) = 2k+1.
\]

Even though the upper bound $\mur(G)+\mur(H)+1$ on $\mur(G\cup H)$ may fail, there is an upper bound for the minimum universal rank of the union of graphs using the minimum universal rank of one and the number of vertices of the other.

\begin{proposition}\label{eincolspace}
For an $n \times n$ symmetric matrix $A$, if $e_n \notin \col(A)$, then $e \in \col(A+\gamma J)$, for all $\gamma\neq 0$.
\end{proposition}
\proof Since $A$ is symmetric, $\mathrm{col}(A)=\mathrm{nul}(A)^{\perp}$. So if $e\notin \mathrm{col}A$, then there exists a vector $x\in \mathrm{nul}(A)$ such that $e^Tx\neq 0$. Thus, $(A+\gamma J)x=\gamma Jx=\gamma (e^Tx)e$, which implies $e \in \mathrm{col}(A+\gamma J)$. \qed

\begin{theorem}
For graphs $G$ and $H$,
\[ 
\mur(G\cup H)\leq \mur(G)+ |V(H)|+1.
\]
\end{theorem}
\proof Assume that $G$ and $H$ have $m$ and $n$ vertices, respectively, and let $\mur(G) =
\rank(U(1,\beta,\gamma,\delta))$. Order the vertices of $G\cup H$ such that the vertices
of $G$ are the first $m$ vertices. Then
\[ 
U_{G\cup H}(1,\beta,\gamma, \delta) = \left[ \begin{array}{c|c}
U_G & \gamma J_{m,n}\\ 
\hline
\gamma J_{n,m} &U_H \end{array}\right].
\] 
If $\gamma=0$, then clearly $\mur(G\cup H)\leq  \mur(G)+ |V(H)|.$
If $\gamma\neq 0$, we consider two cases:

If $e\in \mathrm{col}(U_G)$, then by a similar method used in the proof of Lemma~\ref{lem:upperboundunion}, the matrix $U_{G\cup H}$ can be reduced into the following form
\[ 
 U_1=\left[ \begin{array}{c|c}
U_G & 0\\ 
\hline
0 &U_H+pJ \end{array} \right ]
\] 
for some nonzero number $p$. Therefore,
\[\mur(G\cup H)\leq \mathrm{rank}(U_1)\leq   \mur(G)+ |V(H)|+1.\]

If $e\notin \mathrm{col}(U_G)$, then using Proposition~\ref{eincolspace}, $e \in \mathrm{col}(U_G-\gamma J)$. Let $U'=U_G-\gamma J$, and subtract the $(n+1)$-st column of $U_{G\cup H}(1,\beta,\gamma, \delta)$ from each of the first $n$ columns. The result is the following matrix
\[ 
U_2=\left[ \begin{array}{c|c}
U'& \gamma J_{m,n}\\ 
\hline
R &U_H\end{array}\right].
\]
Since $e \in \mathrm{col}(U')$, the matrix $U_2$ can be reduced to 
\[ 
U_3=\left[ \begin{array}{c|c}
U'&0\\ 
\hline
R &S\end{array}\right].
\] Using the fact that, $\mathrm{rank}(U')\leq \mur(G)+1$, we have 
\[\mur(G\cup H)\leq \mathrm{rank}(U_3)\leq \mathrm{rank}(U')+|V(H)|\leq \mur(G)+|V(H)|+1. \qed\]

\section{Minimum Universal Rank Spread}\label{Sec:murspread}

The {\it mur-spread} of a graph $G$ at vertex $v$, denoted by $\mur_v(G)$, is defined to be
$\mur(G)-\mur(G\setminus \{v\})$. The
following theorem establishes upper and lower bounds for the
mur-spread of a vertex.

\begin{theorem} \label{thm:rankspread}
If a vertex $v$ of $G$ has degree $d$, then
\[-d\leq \mur_v(G) \leq d+2.\]
\end{theorem}
\proof
  Let $U_G=U_G(1,\beta,\delta,\gamma)$ be a universal adjacency matrix associated with the graph $G$.  Let
  $B$ be the submatrix of $U_G$ obtained by deleting the row and
  column of $U_G$ associated with the vertex $v$. Then $U_G$ has the following block form:
\[ U_G(1,\beta,\delta,\gamma)= \left [ 
  \begin{array}{ c|c }
     B & V  \\  \hline
    V^T & {\beta\!+\!\gamma\!+d \delta} \\
  \end{array}
  \right ].
\]
Evidently,
\begin{equation}\label{rankU_Gvs.B}
\rank(U_G)-2 \leq \rank(B) \leq \rank(U_G).
\end{equation}
Let $N(v)$ denote the set of neighbours of $v$ and $D'$ be a diagonal
matrix of the same size as $B$ whose diagonal entry $D'_{ii}$ is $1$
if $v_i\in N(v)$ and $0$ otherwise. Thus there is a universal adjacency matrix for
the graph obtained by removing $v$ from $G$, namely
$U_{G\setminus\{v\}}$, such that $B=U_{G\setminus\{v\}}+\delta D'$.

Using the subadditivity property of rank of matrices, we have the following inequalities
\begin{equation}\label{rankYvs.B}
\rank(B)-d\leq \rank(U_{G\setminus\{v\}}) \leq \rank(B)+d.
\end{equation}

Using equations (\ref{rankU_Gvs.B}) and (\ref{rankYvs.B}), we have
\begin{equation} \label{rankU_Gvs.Y}
\rank(U_G)-(d+2)\leq \rank(U_{G\setminus\{v\}})\leq \rank(U_G)+d.
\end{equation}

Now in (\ref{rankU_Gvs.Y}), if $U_G$ be a universal adjacency matrix with $\rank(U_G)= \mur(G)$, then
\[
\mur(G\setminus \{v\})\leq \rank(U_{G\setminus\{v\}})\leq \mur(G)+d.
\]
which implies $-d\leq \mur_v(G)$. If $U_{G\setminus\{v\}}$ be a universal adjacency matrix with $\rank(U_{G\setminus\{v\}})= \mur(G\setminus \{v\})$ in (\ref{rankU_Gvs.Y}), then
\[
\mur(G)-(d+2)\leq \rank(U_G)-(d+2)\leq  \rank(U_{G\setminus\{v\}})= \mur(G\setminus \{v\}),
\] which implies $ \mur_v(G) \leq d+2.$ \qed

\begin{corollary}
If a vertex $v$ of $G$ has degree $d$, then
\[\max\{-d, -(n-d-1)\} \leq \mur_v(G) \leq \min\{d+2,(n-d-1)+2\}.\]
\end{corollary}
\proof Since 
\[
\mur(G) - \mur(G\setminus \{v\}) = \mur(\overline{G}) - \mur(\overline{G\setminus \{v\}})
          = \mur(\overline{G}) - \mur(\overline{G} \setminus \{v\}), 
\]
simply apply Theorem~\ref{thm:rankspread} to $\overline{G}$, noting
that the degree of $v$ in $\overline{G}$ is $n-d-1$.  \qed

In particular, the following holds:
\begin{corollary}
If  $v$ is a pendant vertex of the graph $G$, then
\[
-1\leq \mur_v(G) \leq 3.\qed
\]
\end{corollary}

\begin{example}
The following examples show that there is a graph with mur-spread $k$
for $-1\leq k \leq 2$. It is an open question to
find a graph with mur-spread equal to $3$ at a pendant vertex.
\begin{enumerate}
\item If $r=s\geq 2$, then $\mur_v(K_r\cup \overline{K_s})\vee
  \{v\}=-1$, for every pendant vertex $v$ (see
  Theorem~\ref{KrcupKsveevspecial} for a proof of this claim).
\item $\mur_v(K_{1,n})=0$ for any pendant vertex $v$ (see Theorem \ref{mur=1}).
\item For $n>2$, $\mur_v(P_n)=1$, for the end-point vertices $v$. 
\item A generalized star is a tree with at most one vertex of degree
  greater than or equal to three. If $G$ is a generalized star on five
  vertices with the degree sequence $1,1,1,2,3$, by calculation it can be shown that $\mur(G)=2$. So if $v$ is a pendant vertex of $G$ whose deletion leaves a star,
using Theorem \ref{mur=1}, we have $\mur_v(G)=2$.
\end{enumerate}
\end{example}

The first example in the list is perhaps the most interesting, since
it is an example of a graph that has a vertex that when removed
leaves a graph that has a strictly larger minimum universal rank. We
will consider this example in more detail in the next section.

\section{Monotonicity}\label{monotonicity}
A parameter for a graph $G$ is called {\it monotone on induced
  subgraphs} if the value of the parameter for the graph is never
smaller than the value on an induced subgraph. In this section, we
show that the minimum universal rank is not in general monotone on
induced subgraphs. To see why this might be true, consider a graph $G$
with a universal adjacency matrix $U$. For any induced subgraph $H$ of $G$,
there is a submatrix of $U$ formed by taking all the rows and columns
corresponding to vertices in $H$. If this submatrix is a universal
matrix for $H$, then the minimum universal rank of $H$ will be no
larger than $\mur(G)$, but this submatrix may not be a universal
matrix for $H$. The problem is with the main diagonal entries, since
these entries are based on the degree of the vertex, and a vertex may have different degrees in different subgraphs. To start, we will give an example of a graph that has
an induced subgraph with a larger minimum universal rank. 

\begin{theorem}\label{KrcupKsveev}
For all nonnegative integers $r$ and $s$, if $s-r+1 \neq 0$, then 
\[ \mur( (K_r \cup \overline{K_s}) \vee \{v\} ) \leq 2.\]
Further,
\begin{enumerate}[(a)] 
\item $\mur( (K_r \cup \overline{K_s}) \vee \{v\} ) = 0$ if and only if ($s=0$), or ($s=1$ and $r=0$);\label{iff1}
\item $\mur( (K_r \cup \overline{K_s}) \vee \{v\} ) = 1$ if and only if
  $s\neq 0$ and $r=0$ or $1$. \label{iff2}
\end{enumerate}
\end{theorem}

\proof
We first prove the statements (\ref{iff1}) and (\ref{iff2}). The graph $(K_r \cup \overline{K_s}) \vee \{v\}$ is
isomorphic to a complete graph if and only if
either $s=0$ or $s=1$ and $r=0$. Since these are the only two
cases in which $(K_r \cup \overline{K_s}) \vee \{v\}$ is a complete
graph, by Theorem~\ref{mur=0}, (\ref{iff1}) holds.

If $r=0$ or $r=1$, then the graph $(K_r \cup \overline{K_s}) \vee
\{v\}$ is a star with $r+s$ edges. Provided that $s>1$, by
Theorem~\ref{mur=1} in the following section, the minimum universal rank of these graphs is one,
thus (\ref{iff2}) holds.

To show the general statement, assume that none of the above cases are satisfied. So either both of
$r$ and $s$ are greater than or equal to $2$ and $s-r+1 \neq 0$ or
$s=1$ and $r\geq 3$. Consider the following universal adjacency matrix:
\[
U = A+ \left(\frac{-1}{r-1}\right) I + \frac{r-1}{s-r+1}J + \frac{1}{r-1}D,
\]
we claim that the rank of this matrix is $2$.

Order the vertices of the graph so that the first $r$ vertices are the
vertices of $K_r$, the next $s$ vertices are the vertices of
$\overline{K_s}$ and the final vertex is $v$. The first $r$ diagonal
entries are all $\frac{s}{s-r+1}$, the next $s$ diagonal entries are
all equal to $\frac{r-1}{s-r+1}$ and the final entry on the diagonal
is 
\[
\frac{s+r-1}{r-1} + \frac{r-1}{s-r+1} = \frac{s^2}{(r-1)(s-r+1)}.
\]
Adjacent vertices have the entry $\frac{s}{s-r+1}$ and nonadjacent
vertices have the entry $\frac{r-1}{s-r+1}$. Thus, the matrix $U$ can
be written as

\[U=\left [
\begin{array}{c|c|c} 
\quad \frac{s}{s-r+1} J_{r,r} \quad  & \quad \frac{r-1}{s-r+1}J_{r,s} \quad &\quad  \frac{s}{s-r+1}e_r\quad \\
 \hline 
\quad \frac{r-1}{s-r+1} J_{s,r} \quad  &\quad  \frac{r-1}{s-r+1}J_{s,s}\quad  &\quad  \frac{s}{s-r+1}e_s\quad \\
 \hline 
\quad \frac{s}{s-r+1}e_r^T \quad  & \quad \frac{s}{s-r+1}e_s^T\quad  & \quad \frac{s^2}{(r-1)(s-r+1)}\quad \\
\end{array}
\right].
\]
Since $s-r+1\neq 0$, the final row is a multiple of the rows in the middle block, which implies $U$ has rank $2$. \qed

\begin{theorem}\label{KrcupKsveevspecial}
For integers $r$ and $s$, if $s \geq 3$ and $s-r+1 = 0$
then 
\[
\mur( (K_r \cup \overline{K_s}) \vee \{v\} ) = 3.
\]
\end{theorem}
\proof As in the proof of Theorem \ref{KrcupKsveev}, order the vertices so
that the first $r$ vertices are from $K_r$, the next $s$
vertices are from $\overline{K_s}$ and $v$ is the last vertex.

Then, any universal adjacency matrix $U=\alpha A+ \beta I +\gamma J +\delta D$ for this graph has
the form
\[ \left [
\begin{array}{c|c|c}
\quad (\gamma\!+\!1) J_{r,r} \!+\! (\beta\!+\!r\delta \!-\! 1)I_{r,r} \quad & \quad \gamma J_{r,s} \quad & \quad (\gamma\!+\!1)e_r \quad  \\
 \hline 
\quad  \gamma J_{s,r}  \quad & \quad \gamma J_{s,s}\!+\! (\beta\!+\!\delta)I_{s,s} \quad & \quad (\gamma\!+\!1)e_s \quad \\
 \hline 
\quad (\gamma\!+\!1) e_r^T \quad & \quad  (\gamma\!+\!1) e_s^T \quad  & \quad \gamma \!+\! \beta\!+\!(r\!+\!s)\delta \quad  \\
\end{array}
\right ]
\]
which can be row reduced to the following matrix
\[
\left [
\begin{array}{c|c|c}
\quad (\beta\!+\!r\delta \!-\! 1)I_{r,r} \quad & \quad (-1) J_{r,s} \quad & \quad [1-(\beta\!+\!(r\!+\!s)\delta)]e_r\quad  \\
 \hline 
\quad  \gamma J_{s,r}  \quad & \quad \gamma J_{s,s}\!+\! (\beta\!+\!\delta)I_{s,s} \quad & \quad (\gamma\!+\!1)e_s \quad \\
 \hline 
\quad (\gamma\!+\!1) e_r^T \quad & \quad  (\gamma\!+\!1) e_s^T \quad  & \quad \gamma \!+\! \beta\!+\!(r\!+\!s)\delta \quad  \\
 \end{array}
\right ].
\]
If $\beta+r\delta - 1 \neq 0$, then the rank of $U$ is at least $r$
which is greater than four. So we assume that $\beta+r\delta - 1
= 0$, and further reduce the matrix to
\[
\left [
\begin{array}{c|c|c}
\quad  0_{r,r} \quad & \quad (-1) J_{r,s} \quad & \quad [1- (\beta\!+\!(r\!+\!s)\delta)]e_r \quad  \\
 \hline 
\quad  \gamma J_{s,r}  \quad & \quad (\beta\!+\!\delta)I_{s,s} \quad & \quad (\gamma\!+\!1)e_s \quad \\
 \hline 
\quad (\gamma\!+\!1) e_r^T \quad & \quad  0_{1,s} \quad  & \quad \gamma \!+\! \beta\!+\!(r\!+\!s)\delta \quad  \\
\end{array}
\right ].
\]
Since $s\geq3$, if $\beta+\delta \neq 0$ then the rank is at least
$3$, so we also assume that $\beta +\delta =0$. With this assumption and the
assumption that $\beta+r\delta - 1 = 0$ we have $\delta =
\frac{1}{r-1}$ and $\beta = -\delta$. The matrix then can be further
reduced to
\[
\left [
\begin{array}{c|c|c}
\quad  0_{r,r} \quad & \quad  (-1) J_{r,s} \quad & \quad -e_r \quad  \\
 \hline 
\quad  \gamma J_{s,r}  \quad & \quad 0_{s,s} \quad & \quad e_s \quad \\
 \hline 
\quad (\gamma\!+\!1) e_r^T\quad & \quad  0_{1,s} \quad  & \quad 1 \quad  \\
\end{array}
\right ].
\]
We also used the facts that $1-(\beta+(r+s)\delta)=-1$, and $\frac{s+r-1}{r-1}=2$. Since there does not exist a value of $\gamma$ such that
$\frac{\gamma+1}{\gamma} =1 $
this matrix has rank $3$.\qed

This particular graph is of interest since it shows that the minimum
universal rank of a graph is not monotone on induced subgraphs. For
example, $G_1=(K_4 \cup \overline{K_3})\vee \{v\}$ is an induced subgraph
of $G_2=(K_4 \cup \overline{K_4})\vee \{v\}$ but $3=\mur(G_1)>\mur(G_2)=2$. (This example also shows that contraction of an edge of a graph can increase the minimum universal rank of a graph.) However, the
minimum universal rank of a graph is monotone under certain conditions. 

\begin{theorem}\label{MonotonicityForRegulars}
  If the minimum universal rank of a graph $G$ is attained with a universal adjacency matrix of $G$ with $\delta=0$, then for any induced subgraph $H$ of $G$ 
\[
\mur(H)\leq \mur(G).
\]
\end{theorem}
\proof
  Let $U=A+\beta I+\gamma J$ be a universal adjacency matrix for $G$
  that attains the minimum rank.  Assume $H$ is obtained from $G$ by
  deleting the set of vertices $R=\{u_1,\ldots,u_r\}$. Then the
  principal submatrix, say $U_H$, of $U_G$ obtained by deleting
  the rows and columns corresponding to $R$, is a universal
  matrix for $H$. Since $\rank(U_H)\leq$ $\rank(U_G)$,
\[
\mur(H)\leq \rank(U_H)\leq \rank(U_G)=\mur(G).\qed
\]

By the discussion in the beginning of the Section~\ref{Sec:regular},
the universal adjacency matrix of a regular graph $G$ can always be written in
the form $U=A+\beta I+\gamma J$. Therefore,
Theorem~\ref{MonotonicityForRegulars} implies the following.

\begin{corollary}
If $G$ is a regular graph and $H$ is an induced subgraph of $G$, then 
\[
\mur(H)\leq \mur(G).\qed
\]
\end{corollary}

This corollary can be used to compute minimum universal rank of some
graphs, for example it can be used to determine the minimum universal
rank of the union of complete graphs with arbitrary sizes. The
complement of such a graph is a complete multipartite graph, so this
will also give the minimum universal rank of these graphs as well.

\begin{theorem}\label{UnionOfCompleteGraphs}
For any integer $k$ and integers $n_1,\ldots,n_k> 1$,
\[
\mur\left(\bigcup_{i=1}^k K_{n_i}\right)=\mur\left(K_{n_1,\ldots,n_k}\right)=k-1.
\]
\end{theorem}

\proof  Let $G=\cup_{i=1}^k K_{n_i}$ and $n=\max\{n_1,\ldots,n_k\}$. Define $G'= \cup_{i=1}^k K_n$, then $G$ is an
  induced subgraph of $G'$. By Theorem~\ref{MonotonicityForRegulars},
  $\mur(G)\leq \mur(G')$.

  The eigenvalues of the adjacency matrix of $G'$ are $n-1$ with
  multiplicity $k$, and $-1$ with multiplicity $k(n-1)$. Therefore,
  using Theorem~\ref{maxmultregular}, we have
\[
\mur(G')=|V(G')|-\left(k(n-1)+1\right)=kn-(kn-k+1)=k-1.
\]
Thus
\begin{align}\label{UpperBound}
\mur(G)\leq k-1.
\end{align}

If we order the vertices of $G$ so that the vertices in $K_{n_i}$ come before the vertices in $K_{n_{i+1}}$, 
then any universal adjacency matrix for $G$ has the form
\[ 
U_G(1,\beta,\gamma,\delta)= \left[\begin{array}{c|c|c|c}
V_1 & \gamma J_{n_1,n_2} & \cdots & \gamma J_{n_1,n_k} \\
 \hline
\gamma J_{n_2,n_1} & V_2 & \cdots & \gamma J_{n_2,n_k} \\
 \hline
\vdots & \vdots & \ddots & \vdots \\
  \hline
\gamma J_{n_k,n_1} & \gamma J_{n_k,n_2} & \cdots & V_k\\
\end{array}\right],
\]
where for any $i=1,\ldots,k$, the $n_i\times n_i$ matrix $V_i$ is as follows:
\medskip
\[
V_i=
\left[
\begin{array}{cccc}
\beta+\gamma+(n_i-1)\delta & \gamma+1& \cdots & \gamma+1  \\
\gamma+1& \beta+\gamma+(n_i-1)\delta & \cdots & \gamma+1  \\
\vdots & \vdots & \ddots & \vdots \\
\gamma+1 & \gamma+1 & \cdots & \beta+\gamma+(n_i-1)\delta \\
\end{array}
\right].
\]
\medskip
Since $n_i>1$, the $k\times k$ submatrix  of $U_G$ that corresponds to the rows
$$\{1,\,\,n_1+1,\,\,n_1+n_2+1,\,\,\ldots,\,\,n_1+\cdots+n_{k-1}+1\},$$
and columns
$$\{2,\,\,n_1+2,\,\,n_1+n_2+2,\,\,\ldots,\,\,n_1+\cdots+n_{k-1}+2\}$$
is
\renewcommand{\arraystretch}{1}
\[
\left[
\begin{array}{cccc}
 \gamma+1& \gamma & \cdots & \gamma  \\
\gamma& \gamma+1 & \cdots & \gamma  \\
\vdots & \vdots & \ddots & \vdots \\
\gamma & \gamma & \cdots & \gamma+1 \\
\end{array}
\right],
\]
subtracting the last column from the previous columns results in the following matrix
\[
\left[
\begin{array}{ccccc}
 1& 0 & \cdots &0& \gamma  \\
0& 1 & \cdots & 0&\gamma  \\
\vdots & \vdots & \ddots & \vdots \\
0 & 0 & \cdots &1& \gamma \\
-1 & -1 & \cdots &-1& \gamma+1 \\
\end{array}
\right],
\]
whose rank is at least $k-1$. This means that $\mur(G)\geq k-1$ and using (\ref{UpperBound}) the result follows. \qed

Note that Theorem~\ref{UnionOfCompleteGraphs} is not true if we drop the condition $n_i>1$.
For example, in the next section we show that $\mur(K_n \cup K_1
\cup K_1)=1$.

\section{Graphs with Minimum Universal Rank Equal to One}\label{Sec:GWSM}

In this section we characterize all graphs $G$ with $\mur(G)=1$.

\begin{theorem}\label{mur=1}
  Let $G$ be a graph with $n=|V(G)|>2$, then $\mur(G)=1$ if and only
  if $G$ or $\overline{G}$ is either $K_r\cup K_s$ for positive $r,s$, with $r+s>2$, or $K_r\cup \overline{K_s}$ for $r,s$ with $1\leq r<n$.
\end{theorem}
\proof
  According to Theorem~\ref{UnionOfCompleteGraphs}, if $G$ or
  $\overline{G}$ is $K_r\cup K_s$ for some $r,s$, with $r+s>2$, then
  $\mur(G)=1$.  Furthermore, if $G=K_r\cup \overline{K_s}$ for some $r,s$ with
  $1<r<n$, then the universal adjacency matrix of $G$ with parameters
  $\alpha=1$, $\beta=0$, $\gamma=0$ and $\delta=\frac{1}{r-1}$ is
  of the form
\[ 
U=A+\frac{1}{r-1}D=\left [
\begin{array}{c|c}
J&0\\\hline
0&0\\
\end{array}
\right ],
\]
whose rank is 1. Thus $\mur(G)\leq 1$, using
Theorem~\ref{mur=0} implies that $\mur(G)=\mur(\overline{G})=1$.

To prove the converse, without loss of generality, assume that $G$ is connected and $\mur(G)=1$. So there are real numbers $\beta,\gamma,\delta$ such that the rank of $U=A+\beta I + \gamma J +\delta D$ is one. Let $S$ be the largest
independent set in $G$. Since $G$ is not a complete graph $s=|S|\geq 2$. Order the vertices of $G$ so that a universal matrix $U=A+\beta I + \gamma J +\delta D$ is of the following form,

\renewcommand{\arraystretch}{1.25}
\[
U= \left [
  \begin{array}{@{}c@{}c@{}c@{}c@{}|@{}c@{}c@{}c@{}c@{}}
    \beta\!+\!\gamma\!+\!\delta d_{v_1} & \gamma & \cdots & \gamma & \ast & \ast &\cdots & \ast \\
    \gamma & \beta\!+\!\gamma\!+\!\delta d_{v_2} & \cdots & \gamma & \ast & \ast& \cdots & \ast \\
    \vdots & \vdots & \ddots & \vdots & \vdots &\vdots & \ddots & \vdots \\
    \gamma & \gamma & \cdots & \beta\!+\!\gamma\!+\!\delta d_{v_s} & \ast & \ast & \cdots & \ast \\ \hline
    \ast & \ast & \cdots & \ast & \beta\!+\!\gamma\!+\!\delta d_{v_{s+1}} & \ast & \cdots & \ast  \\
    \ast & \ast & \cdots & \ast & \ast & \beta\!+\!\gamma\!+\!\delta d_{v_{s+2}} & \cdots & \ast \\
    \vdots & \vdots & \ddots & \vdots & \vdots & \vdots & \ddots & \vdots\\
    \ast & \ast & \cdots & \ast & \ast & \ast & \cdots & \beta\!+\!\gamma\!+\!\delta d_{v_{n}} \\
  \end{array}
  \right ].
\]
 
We claim that all the vertices in $S$ have the same set of
neighbours. To show this, suppose that a vertex $x \in V(G)\setminus
S$ is adjacent to $u\in S$ but not to $v\in S$. As the determinant of

\[U[\{u,v\},\{v,x\}]=\left[ \begin{array}{cc}
\gamma & \gamma+1 \\
\beta+\gamma+\delta d_v & \gamma  \end{array} \right]\]
must be zero, we have
\begin{equation}\label{eq:det0}
\beta(\gamma+1)+\gamma+\delta d_v (\gamma+1)=0.
\end{equation}
Since $G$ is connected, there is a vertex $y\in V(G)\setminus S$
adjacent to $v$. If $y$ is adjacent to $u$, then we have the
following submatrix in $U$:
\[
U[\{u,v\},\{x,y\}]=\left[ \begin{array}{cc}
    \gamma+1 & \gamma+1 \\
    \gamma & \gamma+1 \end{array} \right],
\] 
whose determinant being zero implies that $\gamma=-1$. Substituting
this in (\ref{eq:det0}) leads to a contradiction. Therefore, $y$
cannot be adjacent to $u$. Thus the above submatrix of $U$ is, as
follows
\[U[\{u,v\},\{x,y\}]=\left[ \begin{array}{cc}
\gamma+1 & \gamma \\
\gamma & \gamma+1  \end{array} \right].\]
Since this matrix is singular we have $\gamma=-\frac{1}{2}$. 
Then 
\[U[\{u,v\},\{u,x\}]=\left[ \begin{array}{cc}
\beta+\gamma+\delta d_u & -1/2 \\
-1/2 &1/2  \end{array} \right].\]
Again using singularity we have $\beta+\gamma+\delta d_u=1/2$. Therefore, for any $z\neq u$,  
\[U[\{u,z\},\{u,z\}]=\left[ \begin{array}{cc}
    \frac{1}{2} & \pm\frac{1}{2} \\
    \pm\frac{1}{2} & \beta+\gamma+\delta d_z \end{array} \right],\]
whose singularity results in that fact that $\beta+\gamma+\delta
d_z=\frac{1}{2}$. That is, all diagonal entries are equal to $\frac{1}{2}$.

Assume that $x_1$ and $x_2$ are adjacent vertices of $G$ such that $u$
is adjacent to $x_1$, but $u$ is not adjacent to $x_2$, then $U$
has the following submatrix:
\[
U[\{u,x_1\},\{x_1,x_2\}]=\left[ \begin{array}{cc}
    1/2 & -1/2 \\
    1/2 & 1/2 \end{array} \right],
\] 
which is a contradiction with the rank of $U$ being $1$. Thus, if $u$ is adjacent to a vertex $x_1$, then it must be adjacent to
all the neighbors of $x_1$. But since $G$ is connected, there is a
path $u,x_1,x_2,\cdots,x_t,v$ which implies that $u$ is adjacent to
$v$. This is a contradiction since $u$ and $v$ are both in the
independent set $S$, and so all the vertices in $S$
have the same set of neighbors.

As a result, if a vertex $z\in V(G)\setminus S$ is not adjacent to
a vertex in $S$, then it is not adjacent to any of the vertices in $S$. So $S\cup
\{z\}$ is an independent set, which contradicts the maximality of
$S$. Therefore, all the vertices in $S$ are adjacent to
all the vertices in $V(G)\setminus S$. This implies that $U$ is of the
following form:
\[
{\small
U= \left [
  \begin{array}{@{}c@{}c@{}c@{}c@{}|@{}c@{}c@{}c@{}c@{}}
    \beta\!+\!\gamma\!+\!\delta d_u & \gamma & \cdots & \gamma & \gamma\!+\!1 & \gamma\!+\!1 &\cdots & \gamma\!+\!1  \\
    \gamma & \beta\!+\!\gamma\!+\!\delta d_u & \cdots & \gamma & \gamma\!+\!1 & \gamma\!+\!1& \cdots & \gamma\!+\!1  \\
    \vdots & \vdots & \ddots & \vdots & \vdots &\vdots & \ddots & \vdots \\
    \gamma & \gamma & \cdots & \beta\!+\!\gamma\!+\!\delta d_u & \gamma\!+\!1 & \gamma\!+\!1 & \cdots & \gamma\!+\!1 \\ \hline
    \gamma\!+\!1 & \gamma\!+\!1 & \cdots & \gamma\!+\!1 & \beta\!+\!\gamma\!+\!\delta d_1 & \ast & \cdots & \ast  \\
    \gamma\!+\!1 & \gamma\!+\!1 & \cdots & \gamma\!+\!1 & \ast & \beta\!+\!\gamma\!+\!\delta d_2 & \cdots & \ast  \\
    \vdots & \vdots & \ddots & \vdots & \vdots & \vdots & \ddots & \vdots \\
    \gamma\!+\!1 & \gamma\!+\!1 & \cdots & \gamma\!+\!1 & \ast & \ast & \cdots & \beta\!+\!\gamma\!+\!\delta d_r \\
  \end{array}
  \right ].
  }
\]
Finally, suppose that there is a vertex $x_1\in V(G)\setminus S$ that
is adjacent to a vertex $x_2\in V(G)\setminus S$ but not adjacent to a vertex
$x_3\in V(G)\setminus S$. In this case, the singularity of
\[U[\{u,x_1\},\{x_2,x_3\}]=\left[ \begin{array}{cc}
 \gamma+1& \gamma+1 \\
 \gamma+1& \gamma \end{array} \right]\]
results in $\gamma=-1$. But then we have
\[U[\{u,x_1\},\{v,x_1\}]=\left[ \begin{array}{cc}
 -1& 0 \\
 0& \beta-1+d_{x_1} \end{array} \right],\]
 whose determinant being zero implies that
$$\beta-1+d_{x_1}=0,$$
and so $U$ includes the following submatrix:
\[
U[\{x_1,x_3\},\{x_1,x_3\}]=\left[ \begin{array}{cc}
    0& -1 \\
    -1& \beta-1+d_{x_3} \end{array} \right],
\] 
which has rank two, a contradiction. Hence, the subgraph induced by
$V(G)\setminus S$ is either $K_s$ or $\overline{K_s}$. In the first case
$G=S\vee K_s$ (whose complement is $K_r \cup \overline{K_s}$) and in the second
case $G=K_{r,s}$ (whose complement is $K_r\cup K_s$). \qed

Using Theorem 10 in \cite{HO}, one can provide an alternative method to
prove the ``only if'' part of Theorem~\ref{mur=1}.  Indeed,
$\mur(G)=1$ implies that there exists a universal adjacency matrix for $G$ which
has exactly two distinct eigenvalues; namely 0 and a simple eigenvalue
$\lambda\neq 0$. 

\section{Graphs with large minimum universal rank}

It is known that the only graphs, whose minimum rank is one less than the number of
vertices of the graph are the paths (this is the maximum possible minimum rank). For the
case of minimum universal rank, the maximum possible value is two less than the number of
vertices. It is an interesting question to ask which graphs on $n$ vertices have the maximum minimum
universal rank $n-2$? We have seen that the paths and paths with an isolated vertex
achieve the maximum minimum universal rank; see Example \ref{lem:paths}. Are there any other graphs that also have the maximum
possible minimum universal rank? We consider the paths with an additional edge. 
Define $P'_n$ to be the following graph:
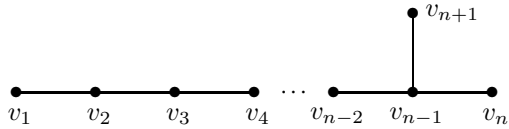
\begin{center}
\unitlength=0.75pt
\begin{picture}(250,50)
\multiput(0,0)(40,0){4}{\circle*{5}}
\multiput(160,0)(40,0){3}{\circle*{5}}
\put(200,40){\circle*{5}}
\put(0,0){\line(1,0){120}}
\multiput(135,0)(5,0){3}{\circle*{1}}
\put(160,0){\line(1,0){80}}
\put(200,0){\line(0,1){40}}
\put(2,-12){\makebox(0,0){$v_1$}}
\put(42,-12){\makebox(0,0){$v_2$}}
\put(82,-12){\makebox(0,0){$v_3$}}
\put(122,-12){\makebox(0,0){$v_4$}}
\put(162,-12){\makebox(0,0){$v_{n-2}$}}
\put(202,-12){\makebox(0,0){$v_{n-1}$}}
\put(242,-12){\makebox(0,0){$v_{n}$}}
\put(220,40){\makebox(0,0){$v_{n+1}$}}
\end{picture}\unitlength=1pt
\bigskip
\medskip
\end{center}

For $n=4,5$ we know that $\mur(P'_n)=n-1$ which is the number of vertices of the graph minus two in each case. So, there are graphs other than paths, with the maximum possible minimum universal rank. But $n=4,5$ are the only cases known for this family of graphs. Indeed, for infinitely many values of $n$ the minimum universal rank of $P'_n$ is three less than the number of vertices. 
\begin{proposition}

\begin{enumerate}[(a)]
\item For $n\geq 3$, if $n\equiv 0\,\,\,$ (mod 3), then $\mur(P'_n)=n-2$.
\item For $n\geq 6$, if $n+1=4k$, then  $\mur(P'_n)=n-2$.
\end{enumerate}
\end{proposition}
\proof Since $P'_n$ has an induced path on $n$ vertices, using Lemma \ref{lbd-indpath}, we have $\mur(P'_n)\geq n-2$. Under the assumption of part (a), the universal matrix $U(1,1,  -\frac{1}{n+1}, -1)$ has rank $n-2$, and under the assumption of part (b), the universal matrix $U(1,0,0,-\frac{1}{2k})$ has rank $n-2$. \qed

Moreover, for $n=8$, the $9\times 9$ universal matrix with parameters  $\alpha=1,-\beta=\delta=\frac{1\pm \sqrt{5}}{2}$ and $\gamma=\frac{1}{3\delta-5}$ has rank $6$. And, for $n=10$, any set of the parameters $\alpha=1,-\beta=\delta=-2\cos \frac{2\pi}{7},\gamma=\frac{-1}{\delta^2-5\delta+6}$ or $\alpha=1,-\beta=\delta=-2\cos \frac{6\pi}{7},\gamma=\frac{-1}{\delta^2-5\delta+6}$ gives the minimum universal rank equal to $8$. This leads us to speculate that $\mur(P'_n)=n-2$, for $n\geq 6$. 

\section*{Acknowledgements}
The work in this paper was a joint project of the Discrete Mathematics Research Group at the University of Regina, attended by all the authors.  They would like to thank the other members of the group (J.~C.~Fisher and A.~M.~Purdy) for their participation.  S.~Fallat and K.~Meagher acknowledge support from an NSERC Discovery Grant. S.~Nasserasr is a PIMS Postdoctoral Fellow. Y. Fan 's research was supported by National Natural Science Foundation of China (11071002), Key Project of Chinese Ministry of Education (210091), Specialized Research Fund
for the Doctoral Program of Higher Education (20103401110002), Anhui Provincial Natural Science Foundation
(10040606Y33), Project of Educational Department of Anhui Province (2009SQRZ017ZD, KJ2010B136),
Academic Innovation Team of Anhui University Project (KJTD001B).


\begin{thebibliography}{10}

\bibitem{BH}
\newblock{E. Brouwer, W. H. Haemers, Spectra of graphs, ebook.}

\bibitem{BR}
\newblock{R. Brauldi and H. Ryser. Combinatorial Matrix Theory. Cambridge University Press, Cambridge, 1991.}

\bibitem{CDS}
\newblock{Cvetkovi\'c, D. M., Doob, M., and Sachs, H. Spectra of Graphs: Theory and Applications, 3rd rev. enl. ed. New York: Wiley, 1998.}

\bibitem{FH}
\newblock{S. Fallat and L. Hogben. The minimum rank of symmetric matrices described by a graph:a survey. Linear Algebra and its Applications, 426:558-582, 2007.}

\bibitem{DHK}
\newblock{E. R. van Dam, W. H. Haemers and J. H. Koolen, Cospectral graphs and the generalized adjacency matrix, Linear Algebra Appl. 423 (2007), 33-41.}

\bibitem{HO}
\newblock {Haemers, Willem H. and Omidi, G.R., Universal Adjacency Matrices with Two Eigenvalues (November 3, 2010).}
\newblock  {CentER Discussion Paper Series No. 2010-119. Available at SSRN: http://ssrn.com/abstract=1717756.}

\bibitem{JN}
\newblock{C. R. Johnson \& M. Newman, A note on cospectral graphs, J. Combin. Th. (B) 28 (1980) 96-103.}

\bibitem{M}
\newblock{R. Merris, Laplacian matrices of graphs; a survey,  Linear Algebra Appl. 197\&198 (1994), 143-176.}

\end{thebibliography}
\end{document}